\documentclass{svproc} 
   \usepackage{url}

\usepackage{mathptmx}   
\usepackage{helvet} 
\usepackage{courier} 
\usepackage{amsmath,amsfonts}  
\usepackage{graphicx}
\usepackage[bottom]{footmisc}           
    
\usepackage{mathptmx}         
\usepackage{helvet}            
\usepackage{courier}            
\usepackage{type1cm}             
\usepackage{makeidx}             
\usepackage{graphicx}          
\usepackage{multicol}          
\usepackage[bottom]{footmisc}  
     
\usepackage{epsfig}     
\usepackage{psfrag,rotating}     
\usepackage{amssymb,floatflt,enumerate} 
\usepackage{amsmath,amscd,psfrag,leqno} 
\usepackage[mathscr]{eucal}  
\usepackage{shadethm}           
\usepackage{graphicx,subfigure}   
\usepackage{overpic,contour}

\contourlength{0.3mm}

\def\phi{\varphi}

\def\RR{{\mathbb R}}

\def\CC{{\mathbb C}}

\def\Vkt#1{{\mathbf #1}}

\definecolor{blau}{rgb}{0,0,1}

\definecolor{rot}{rgb}{1,0,0}

\definecolor{green}{rgb}{0,1,0} 

\definecolor{cyan}{rgb}{0,1,1}

\definecolor{violet}{rgb}{1,.54,1}

\definecolor{magenta}{rgb}{1,0,1}

\begin{document}
\mainmatter              
\title{
On flexes associated with higher-order flexible \\ bar-joint frameworks}
\titlerunning{On flexes associated with higher-order flexible bar-joint frameworks}  
%
\author{Georg Nawratil}
\authorrunning{Georg Nawratil} 
%
%
\institute{Institute of Discrete Mathematics and Geometry \& \\ Center for Geometry and Computational Design, TU Wien, Austria,\\
\email{nawratil@geometrie.tuwien.ac.at},\\ WWW home page:
\texttt{https://www.dmg.tuwien.ac.at/nawratil/}
}

\maketitle              

\begin{abstract}
The famous example of the double-Watt mechanism given by Connelly and Servatius raises some problems concerning the classical definitions of higher-order flexibility and rigidity, respectively. 
Recently, the author was able to give a proper redefinition of the flexion/rigidity order  
for bar-joint frameworks, but the question for the flexes associated with higher-order flexible structures remained open. In this paper we properly define these flexes based on the theory of algebraic curves and demonstrate their computation by means of Puiseux series. The presented algebraic approach also allows to take reality issues into account.
\end{abstract}

\section{Introduction}
\label{sec:introduction}

A bar-joint framework $G(\mathcal{K})$ consists of a knot set $\mathcal{K}=\left\{X_{1},\ldots, X_w\right\}$
and a graph $G$  on $\mathcal{K}$ encoding the combinatorial structure. 
A knot $X_i$ corresponds to 
a rotational/spherical joint (without clearance) in the case of a planar/spatial framework. 
An edge connecting two knots corresponds to a bar. We denote the number of edges by $e$. Note that $w$ and $e$ are finite numbers as we do not consider infinite frameworks.

By defining the lengths of the bars, which are assumed to be non-zero,    
the intrinsic metric of the framework is fixed. 
In general the assignment of the intrinsic metric does not uniquely determine the embedding of the framework into the Euclidean space $\RR^s$, thus such a framework can have 
different incongruent realizations. 

The relation that two elements of the knot set are edge-connected can also be expressed algebraically, which implies $e$ quadratic equations $c_1,\ldots, c_e$ in $m$ unknowns. 
Note for the planar case ($s=2$) we get $m=2w-3$ and for the spatial one ($s=3$) $m=3w-6$, after eliminating the isometries of the complete framework.

If the algebraic variety  $V(c_1,\ldots,c_e)$ is positive-dimensional then the framework is flexible; otherwise rigid.  The framework is called isostatic (minimally rigid) if the removal of any algebraic constraint $c_i$ will make the framework flexible. 
In this case $m=e$ has to hold.

If $V(c_1,\ldots ,c_l)$ is zero-dimensional, then each real solution corresponds to a realization $G(\Vkt X)$ of the framework with 
$\Vkt X=(\Vkt x_1,\ldots , \Vkt x_w)\in\RR^m$. 


\subsection{Review and outline}

The classical definition of a higher-order flex can be written as follows according to Stachel \cite{stachel_proposal}: 

\begin{definition} \label{def1}
A framework realization has a $n^{th}$-order flex if for each vertex $\Vkt x_i$ ($i=1,\ldots,w$) 
there is a polynomial function 
\begin{equation}\label{eq:flex}
\Vkt x_i':=\Vkt x_i+ \Vkt x_{i,1}t+\ldots + \Vkt x_{i,n}t^n \quad \text{with} \quad n>0
\end{equation}
such that
\begin{enumerate}
\item
the replacement of $\Vkt x_i$ by $\Vkt x_i'$ in the equations $c_1,\ldots ,c_e$ gives stationary values 
of multiplicity $\geq n+1$ at $t=0$;
\item
the velocity vectors $\Vkt x_{1,1},\ldots ,\Vkt x_{w,1}$ do not originate from a rigid body motion (incl.\ standstill) 
of the complete framework; i.e.\ they are said to be non-trivial. 
\end{enumerate}
\end{definition}

Of special interest is the flex of maximal order $n^*$ as then the realization under consideration is said to have flexion order $n^*$ and a rigidity of order $n^*+1$ (cf.\ \cite{servatius}), respectively. The famous example of the double-Watt mechanism given by Connelly and Servatius \cite{servatius} raises some problems concerning these definitions, as they attest the cusp configuration a third-order rigidity, which conflicts with its continuous flexion. 

Two attempts are known to the author to resolve this dilemma; namely the one of Gaspar and Tarnai \cite{gaspar} and the one of Stachel \cite{stachel_proposal}. In the following we only discuss the latter approach as 
it can be shown \cite[Rmk.\ 3]{MMT} that it also 
contains the one of Gaspar and Tarnai.
Stachel \cite{stachel_proposal}  tried to settle the problem by following a more general notation of $(k,n)$-flexes suggested by Sabitov \cite{sabitov} which replaces Eq.\ (\ref{eq:flex}) 
by
\begin{equation}\label{eq:flexk}
\Vkt x_i':=\Vkt x_i+ \Vkt x_{i,k}t^k+\ldots + \Vkt x_{i,n}t^n \quad \text{with} \quad n\geq k>0
\end{equation}
where $\Vkt x_{1,k},\ldots ,\Vkt x_{w,k}$ is non-trivial. In addition Stachel assumed that Eq.\ (\ref{eq:flexk}) represents 
an {\it irreducible flex}; this means that  Eq.\ (\ref{eq:flexk}) does not result from 
a polynomial parameter substitution of the form
    $t=\overline{t}^q(a_0+a_1\overline{t}+a_2\overline{t}^2+\ldots)$ 
with $a_0\neq 0$ and $q>1$ into a lower-order flex.
This approach implies for the cusp configuration of the double-Watt mechanism an irreducible $(2,\infty)$-flexibility (see also \cite[ Ex.\ 2]{MMT} for details). 

Stachel did not publish this approach due to a new dilemma arising from the presentation \cite{stachel_aim} of another double-Watt mechanism, 
which is extended by a Kempe-mechanism. 
For the resulting rigid framework Stachel was not able to determine a unique flexion order $n^*$ with his proposed definition as he ended up with an infinite sequence of irreducible $(k,3k-1)$-flexes  
for $k=1,2,\ldots$ (see also \cite[Ex.\ 3]{MMT}). 

Recently, we were able to give a proper redefinition of the flexion/rigidity order 
for bar-joint frameworks  \cite[Defs.\ 3 and 4]{MMT} inspired by Sabitov's finite algorithm for testing the bendability of a polyhedron \cite{sabitov}. For a configuration that does not belong to a continuous flexion of the framework its flexion order $r$ is defined as the number of coinciding framework realizations minus 1; its rigidity order is then defined as  $r+1$. 
But the question for the flexes associated with structures having a flexion order $r>1$ remained open.


In Section \ref{sec:defcomp} we present a modification of Stachel's detailed definition of the $(k,n)$-flexes introduced by Sabitov, which settles the problem to identify the correct value for $k$. This approach is motivated in Section \ref{sec:removal} by the so-called removal procedure for isostatic frameworks. Moreover, in Section \ref{sec:defcomp} we also discuss the computation of $(k,n)$-flexes by means of Puiseux series, which are recapped within the next section, where some preliminary considerations and fundamentals are given.


\subsection{Preliminary considerations and fundamentals}\label{sec:pcf}

The problem in Stachel's approach for computing the $(k,n)$-flexes is that the correct value for $k$ is unknown. Let us shade some light on the meaning of $k$ by considering first the limit case $n=\infty$  ($\Leftrightarrow$ the framework is continuous flexible) as in this case the flexes correspond to branches of an 
algebraic curve (e.g.\ \cite{burau,semple,walker}), 
which can be locally parametrized by Puiseux series. This theory is well established for planar algebraic curves,  
but it can 
also be extended to the non-planar case (e.g.\ \cite{alonso,jensen,maurer}). 
For the approach taken in this conceptual paper it is sufficient to understand the Puiseux series expansion of planar algebraic curves, which is summarized next: 

Without loss of generality we can always assume that the algebraic curve $C$ is moved in a way that 
the considered point of $C$ equals the origin of the reference frame. 
According to the Newton-Puiseux Theorem each branch of a planar algebraic curve $C$, which is given by the zero-set of a polynomial $P(x,y)$, 
can be parametrized by a pair of convergent power series. 
This results from the Theorem of Puiseux that  $P(x,y)=\alpha_{0}(x)+\alpha_{1}(x)y+\ldots +\alpha_{n}(x)y^{n}$ 
can be written as 
$\prod_{u=1}^{n}(y-S_{u}(x))$ where $S_{u}(x)$ is the Puiseux series of the form 
\begin{equation}
S_{u}(x)= \beta_{u,1}x^{\tfrac{\nu_{u,1}}{\nu_{u,0}}}+\beta_{u,2}x^{\tfrac{\nu_{u,2}}{\nu_{u,0}}}+\ldots \quad \text{with}\quad 0<\nu_{u,0},\quad 0<\nu_{u,1}<\nu_{u,2}<\ldots
\end{equation}
and relative prime $\nu_{u,0},\nu_{u,1},\nu_{u,2},\ldots$, which 
can be obtained my means of Newton diagrams.
From this we get the minimal parametrization 
of a branch\footnote{\label{fn1}Note the technical detail that 
this parametrization is only possible for branches which are not contained within the y-axis given by $x=0$.} at the origin of 
 $C$ as
\begin{equation}
x_{u}(t)=t^{\nu_{u,0}}, \quad y_{u}(t)= \beta_{u,1}t^{\nu_{u,1}}+ \beta_{u,2}t^{\nu_{u,2}}+\ldots .
\end{equation} 
This parametrization is not unique for $\nu_{u,0}>1$ as then the substitution of $t$ by $\epsilon t$ yields and 
equivalent parametrization; where $\epsilon$ denotes a $\nu_{u,0}$-th root of unity. 
Therefore the number of branches is less or equal to $n$.
Note that the order of a branch at the origin of $C$ is  $\min(\nu_{u,0},\nu_{u,1})$. 
It is well-known that this order remains invariant under the choice of the coordinate frame, thus an 
ansatz for the branch of order $k$ has to equal Eq.\ (\ref{eq:flexk}). 
This shows the algebraic geometric meaning of $k$ in Stachel's approach for  $n=\infty$.


\section{Removal procedure for isostatic bar-joint frameworks}\label{sec:removal}

From the standpoint of kinematics it is quite natural to remove the $i$th bar of an isostatic bar-joint framework for $i\in\left\{1,\ldots, e\right\}$
and to consider the resulting one degree of freedom mechanism.  Then one can compute in the configuration $\Vkt X$ of interest  
the branches of the one-dimensional configuration curve generated by the ideal 
\begin{equation}\label{ideal:basic}
\langle c_1,\ldots,c_{i-1},c_{i+1},\ldots,c_e\rangle. 
\end{equation}
Then we check, up to which order each branch is compatible with the removed condition $c_i=0$, which contains $\Vkt X$ as regular point,  by determining the intersection multiplicity $n+1$. 
Then this branch implies a $n^{th}$-order flex. 

\begin{remark}
A result on this removal procedure for the ordinary case of $(1,n)$-flexes is given by Stachel in \cite[Lem.\ 1]{stachel_planar}. 
\hfill $\diamond$
\end{remark}

\begin{figure}[t]
\begin{center}
\begin{overpic}
    [height=50mm]{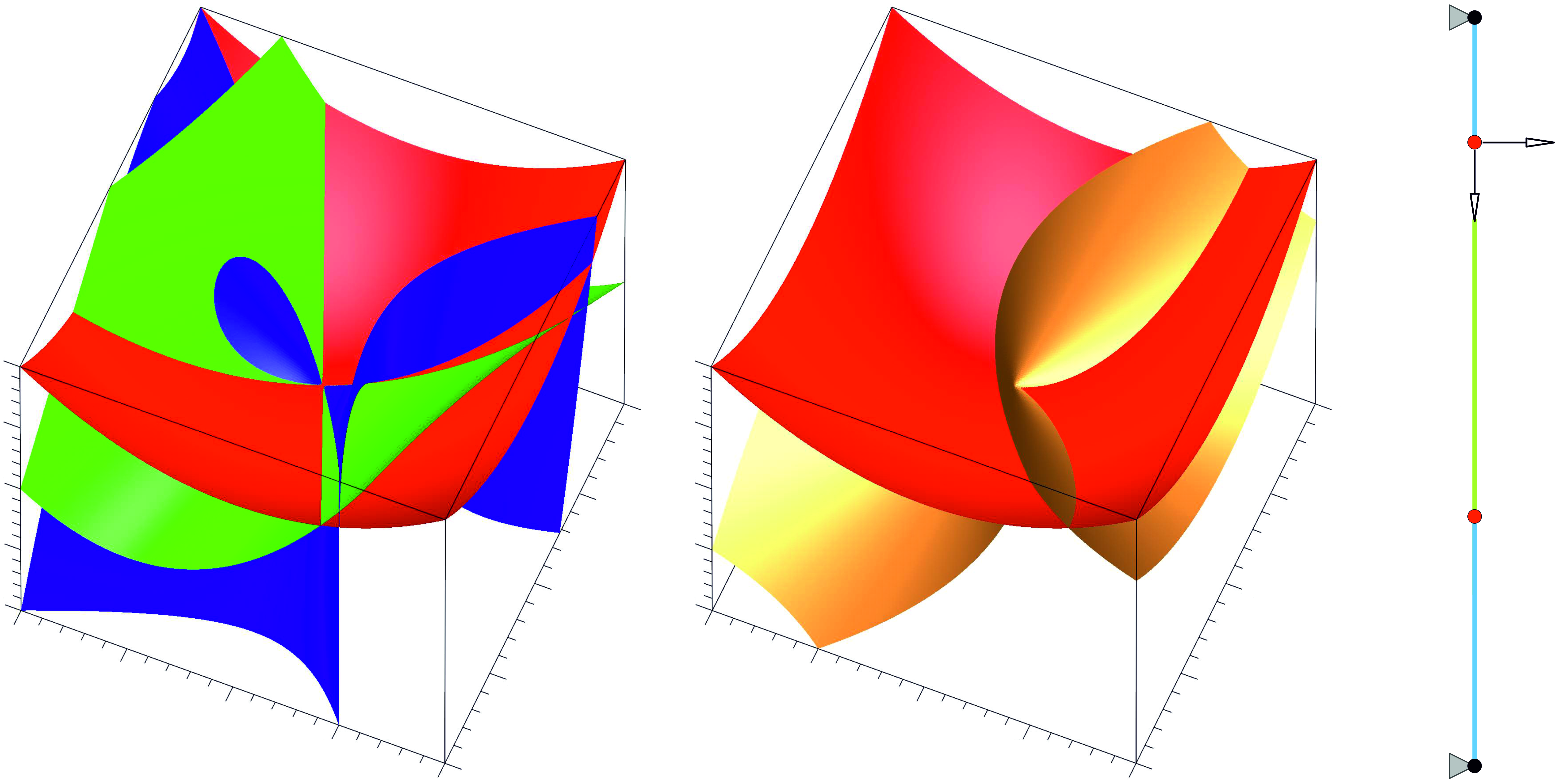}
\begin{scriptsize}
\put(21,0.5){$y$}
\put(13,3){$0$}
\put(33,5){$x$}
\put(36,11.5){$0$}
\put(-1.7,19){$0$}
\put(-1.2,24){$z$}
\put(42.7,18.5){$0$}
\put(43.2,24){$z$}

\put(65.5,0.5){$y$}
\put(57.5,3){$0$}
\put(77.5,5){$x$}
\put(80.5,11.5){$0$}
\put(95.5,48.3){$F_1$}
\put(95.5,0){$F_2$}
\put(99,39.5){$y$}
\put(95.5,36.8){$x$}
\put(91.1,40){$M_1$}
\put(91.1,16){$M_2$}
\put(96,8){$\overline{M_2F_2}=2$}
\put(96,28){$\overline{M_1M_2}=3$}
\put(96,44){$\overline{F_1M_1}=1$}
\end{scriptsize}     
  \end{overpic} 
\end{center}	
\caption{Left: Intersection of the three quadrics $C_1$ (red), $C_2$ (green) and  $C_3$ (blue). Center:   
The intersection curve of the cone $C_4$ (yellow) and $C_1$. 
Right: Immobile 4-bar mechanism.}
  \label{fig1}
\end{figure}
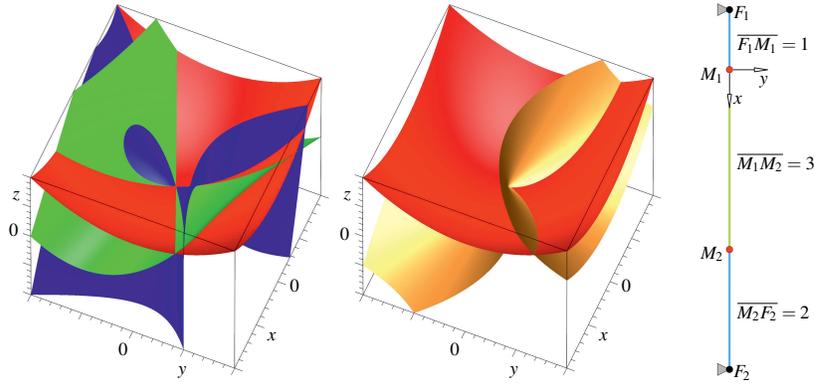

\begin{example}\label{ex2}
For purpose of illustration we do not use an example which corresponds to an actual framework. 
The example is constructed based on the data provided in \cite[Table 1]{tu}.
Let us consider the three quadrics $C_i$ in $\RR^3$ given by $c_i(x,y,z)$ with:
\begin{equation}
\begin{split}
c_1(x,y,z):=&x^2+y^2-2z, \quad
c_2(x,y,z):=y^2+xy-z, \\
c_3(x,y,z):=&2x^2-3xy-2y^2-2yz+z.
\end{split}
\end{equation}
Note that the three regular quadrics intersect in the origin with multiplicity 4 ($\Rightarrow$ $r=3$), which can for example be checked with the method given in \cite[Sec.\ 3.1]{MMT}. 
It can easily be seen that these quadrics have a common tangent plane $z=0$ at the origin, which is also a double point of each three possible intersection curves (cf.\ Fig.\ \ref{fig1}-left). 
By applying only the removal procedure one would end up with the conclusion that there are always two linear branches, which both imply flexes of order $1$. \hfill $\diamond$
\end{example}

But by this kinematic motivated removal procedure no complete picture of the flexes 
can be obtained, which becomes clear by 
looking at the problem from the standpoint of algebraic geometry. 
The ideal of Eq.\ (\ref{ideal:basic}) can be generalized to
\begin{equation}\label{ideal:general}
\langle c_1+\lambda_1c_i,\ldots, c_{i-1}+\lambda_{i-1}c_{i},c_{i+1}+\lambda_{i+1}c_{i},\ldots , c_e+\lambda_ec_i \rangle
\end{equation}
where the $\lambda_1,\ldots,\lambda_{i-1},\lambda_{i+1},\ldots,\lambda_e\in\RR$ imply a $(e-1)$-parametric set of curves, 
whose branches have to be intersected with the  hypersurface $c_i=0$ containing the origin as regular point. The resulting intersection multiplicity $n+1$ yields again the order $n$ of the flex implied by the corresponding branch. 
 Note that if we add $c_i=0$ to the 
ideal  
 of Eq.\ (\ref{ideal:general}), then the resulting ideal 
is equivalent to the 
initial one $\langle c_1,\ldots,c_e\rangle $ 
including also its intersection multiplicity.

\begin{example}\label{ex:con} Let us continue with Example \ref{ex2}.  
Within the bundle of quadrics spanned by $C_1,C_2,C_3$ there also exists a pencil of cones $\mathcal{P}$ having their apexes in the origin (cf.\ \cite{li}). 
Moreover $\mathcal{P}$ touches the plane 
$z=0$ in a pencil of lines through the origin. 
Therefore each cone of $\mathcal{P}$ intersect $C_1$ in a quartic curve having a cusp in the origin in direction of the cone's generator contained in $z=0$. 
This is illustrated in Fig.\ \ref{fig1}-center for the cone $C_4\in\mathcal{P}$ given by $c_4=0$ with $c_4(x,y,z):=x^2-2yz$. 
Therefore there exists also a pencil of $(2,3)$-flexes. 
Note that higher-order cusps with $k>2$ are not possible due to degree reasons\footnote{For a non-planar branch two points of this branch would span with the cusp a plane intersecting the branch (which is maximal of order 4) in at least five points (counted with multiplicity) if $k>2$; a contradiction. A planar branch has to be a conic, which cannot have cusps.}. 
The geometric reasoning given within this example  will also be verified by the detailed algebraic analysis presented in Example \ref{ex3} of the next section.
\hfill $\diamond$
\end{example}


\section{Definition and computation of (k,n)-flexes}\label{sec:defcomp}

The above given considerations for isostatic bar-joint frameworks motivate the following redefinition of a $(k,n)$-flex for a general 
bar-joint framework:

\begin{definition} \label{def_flex}
A bar-joint framework, which is not continuous flexible, has a 1-parametric $(k,n)$-flex if for each vertex $\Vkt x_i$ ($i=1,\ldots,w$) 
there is a polynomial function given in Eq.\ (\ref{eq:flexk})
such that
\begin{enumerate}
\item
the replacement of $\Vkt x_i$ by $\Vkt x_i'$ in the equations $c_1,\ldots,c_e$ of the edge lengths gives stationary values 
of multiplicity $\geq n+1$ at $t=0$;
\item
the vectors $\Vkt x_{1,k},\ldots ,\Vkt x_{w,k}$ are non-trivial; 
\item
Eq.\ (\ref{eq:flexk}) can be extended to a minimal parametrization of a branch of order $k$ of an algebraic curve, 
which corresponds to a one-dimensional irreducible component of a variety determined by an ideal, whose
generators  
are contained in the linear 
family of quadrics spanned by $c_1,\ldots, c_e$.
\end{enumerate}
\end{definition}
Note that Sabitov's original definition \cite[page 230]{sabitov} only includes the items 1 and 2, which was extended by Stachel by the condition that the flex has to be irreducible. 
In Def.\ \ref{def_flex} this additional condition is replaced by item 3 in order to resolve the dilemma of  ending up with an infinite series of possible 
$(k,n)$-flexes, as pointed out by Stachel's example of the extended double-Watt mechanism. 
This is due to the fact that $k$ cannot be greater than\footnote{$r$ is only an upper bound for $k$ as a high flexion order $r$ can also result from multiple branchings into associated flexes of lower order.} the flexion order $r$. 
Note that item 3 of  Def.\ \ref{def_flex} implies a 
global construction for the determination of $(k,n)$-flexes (like the computation of the flexion order $r$ \cite{MMT}).

We added also the assumption that the framework 
is not continuous flexible, 
but Def.\ \ref{def_flex}  also holds for frameworks with a 1-dimensional mobility ($\Leftrightarrow$ $n=\infty$; cf.\ Example \ref{ex4}). 
Moreover, we speak more precisely of a ``1-parametric flex'' as one can also think of $p$-parametric flexes with $p>1$, which are not within the scope of this paper. Their definition can be done similarly to Def.\  \ref{def_flex} and their study can be based on the local approximations of $p$-surfaces by multivariate Puiseux series expansions \cite{aroca,buchacher,mc1,mc2}.

In the remainder of this section we tackle the computation of $(k,n)$-flexes according to Def.\ \ref{def_flex}. 
To do so, we start with the construction of the minimal parametrizations of the branches of the algebraic curves mentioned in item 3 of Def.\ \ref{def_flex}. 
For planar curves this can be done according to Section \ref{sec:pcf} and for 
space curve we use the method given by Melanova \cite[page 112]{melanova}. Her  strategy is to  project a space curve to all possible coordinate planes containing a selected coordinate axis\footnote{Footnote \ref{fn1} implies that the hyperplane orthogonal to this axis through the origin is not allowed to contain a branch of the space curve.}
and to construct the minimal parametrizations of the planar projections. Then  from all these parametrizations a minimal one of the space curve can be obtained (for details see \cite{melanova}). 
In the following this strategy is demonstrated for the example already discussed in Section \ref{sec:removal}, where we use the resultant method for the projection to the coordinate planes.

\begin{example}\label{ex3}
First we look at the ideal $I_1=\langle c_1+\lambda_1c_3,c_2+\lambda_2c_3\rangle$. For $\lambda_1,\lambda_2\neq 0$ this ideal is equivalent with the ideals
$\langle c_1+\mu_1c_2,c_3+\mu_3c_2\rangle$ with $\mu_1,\mu_3\neq 0$
and 
$\langle c_2+\nu_2c_1,c_3+\nu_3c_1\rangle$
with $\nu_2,\nu_3\neq 0$, respectively. 
It can easily be seen that we only have to discuss this ideal $I_1$ as well as the ideals 
    $I_2=\langle c_1+\mu_1c_2,c_3\rangle$ and 
    $I_3=\langle c_2,c_3\rangle$
in order to cover all algebraic curves which can be obtained by the intersection of two quadrics belonging to the bundle spanned by $C_1,C_2,C_3$. 

In order to demonstrate the procedure we start with the discussion of the ideal $I_1$:
First we eliminate $z$ from the two generators of $I_1$ by the resultant method, which yields the polynomial $P(x,y)$ given by:
\begin{equation*}
x^2 - 2xy - y^2 
+ \lambda_1(2x^2 - 2xy - y^2 - 2xy^2 - 2y^3) + 
\lambda_2 (2x^2y + 2y^3  - 5x^2 + 6xy + 3y^2). 
\end{equation*}
By using {\sc Maple} one can compute the Puiseux series in dependence of $\lambda_1,\lambda_2$ under the assumption 
\begin{equation}\label{eq:assume}
 (\lambda_1-3\lambda_2+1)(3\lambda_1-8\lambda_2+2)\neq 0.   
\end{equation}
 In this case we get two linear branches with minimal parametrization
\begin{equation}
x(t)=t, \quad y_{\mp}(t)= -\tfrac{ \lambda_1 - 3\lambda_2 +1 \mp \sqrt{(3\lambda_1 - 8\lambda_2 + 2)(\lambda_1 - 3\lambda_2 + 1)} }{\lambda_1 - 3\lambda_2 + 1}t
+\ldots .
\end{equation} 
Now we eliminate $y$ by the resultant method from the generators of $I_1$ yielding the polynomial $P(x,z)$ with:
\begin{equation}
\begin{split}
 &2x^4 - 4x^2z + z^2
  +\lambda_1\lambda_2(32x^2z + 26xz^2 + 12z^3 -25x^4 - 12x^3z - 4x^2z^2  - 6z^2) +\\
 &\lambda_1( 7x^4+ 2zx^3-10 z x^2-6 z^2 x+2 z^2)
 +2\lambda_2(13x^2z + 4xz^2 -7x^4 - 2x^3z  - 3z^2).
\end{split}
\end{equation}
Again under the assumption of Eq.\ (\ref{eq:assume}) we can compute the  Puiseux series in dependence of $\lambda_1,\lambda_2$. We get again two linear branches with minimal parametrization
\begin{equation}
x(t)=t, \quad z_{\mp}(t)= \tfrac{3\lambda_1 - 7\lambda_2 + 2 \mp \sqrt{(3\lambda_1 - 8\lambda_2 + 2)(\lambda_1 - 3\lambda_2 + 1)} }{\lambda_1 - 3\lambda_2 + 1}t
+\ldots .
\end{equation} 
Then the two linear branches of the space curve have the minimal parametrization $(x(t),y_\mp(t),z_\mp(t))$. Plugging both  parametrizations into $c_3$ yield
\begin{equation*}
 \tfrac{1}{\lambda_1 - 3\lambda_2 + 1}t^2 \mp 
\tfrac{(8\lambda_1 - 20\lambda_2 + 6)\sqrt{(3\lambda_1 - 8\lambda_2 + 2)(\lambda_1 - 3\lambda_2 + 1)} + 2(\lambda_1 - 3\lambda_2 + 1)(4 + 6\lambda_1 - 15\lambda_2)}{\lambda_1 - 3\lambda_2 + 1}^3 t^3 +o(t^3).
\end{equation*}
This shows that each of the two linear branches imply a $(1,1)$-flex. 
For the discussion of the special cases excluded by  Eq.\ (\ref{eq:assume}) we refer to Appendix A, where also the ideals $I_2$ and $I_3$ are discussed in detail. In summary, these computations prove that only $(1,1)$ and $(2,3)$-flexes are possible verifying the argumentation given in Example \ref{ex:con}.
\hfill $\diamond$
\end{example}

The algebraic approach for the computation of the $(k,n)$-flexes  operates over $\CC$ but it also allows to take reality issues into account by using only the real part of the minimal parametrizations of the branches. This can be used to determine the highest real flex $(k_{\max},n_{\max})$, which is of interest as it  complements the flexion order $r$ according to  \cite[item 3 of Sec.\ 5]{MMT}. In the following this is demonstrated for an immobile 4-bar linkage.

\begin{example}\label{ex4}
We use the same dimensioning for the immobile 4-bar mechanism (cf.\ Fig.\ \ref{fig1}-right) as in \cite[Sec.\ 3.2]{nayak} and coordinatize the fixed joints by $F_1=(-1,0)$, $F_2=(5,0)$ and the moving ones by $M_1=(a,b)$, $M_2=(3+c,d)$. 
Moreover, we have the constraints:
\begin{equation*}
c_1:\,\,\|M_1-F_1\|^2-1^2=0,\quad
c_2:\,\,\|M_2-M_1\|^2-3^2=0,\quad
c_3:\,\,\|M_2-F_2\|^2-2^2=0.
\end{equation*}
The parametrization of the space curve splits up into two conjugate complex linear branches (see Appendix B for the detailed computation) given by $b(t)=t$ and
\begin{equation*}
a(t)=\tfrac{-1}{2}t^2  -\tfrac{1}{8}t^4  + \ldots, \,\,\,\,\,\,\,\,
c(t)=\tfrac{-8 \pm 6I}{25}t^2  + \tfrac{-79 \pm 3I}{1250}t^4 + \ldots,  \,\,\,\,\,\,\,\,
d(t)=\tfrac{2\pm 6I}{5}t +  \tfrac{6\pm 33I}{250}t^3 +\ldots
\end{equation*}
through the origin, which is the only real point. 
By plugging the real part of this parametrization 
into $c_1,c_2,c_3$ it can be seen that  at least $t^2$ factors out from the resulting three expressions. This implies that the triple $(r;k_{\max},n_{\max})$ equals $(\infty;1,1)$. \hfill $\diamond$ 
\end{example}

Finally, it should be noted that for the central example of the paper (cf.\ Examples \ref{ex2}--\ref{ex3}) the triple $(r;k_{\max},n_{\max})$ equals $(3;2,3)$. In contrast the original double-Watt mechanism has the triple $(\infty;2,\infty)$ according to the analysis given in \cite{servatius,gaspar,muller,MMT,nayak,stachel_proposal}.

\section{Conclusion and future work}

In Def.\ \ref{def_flex}  we gave a modification of Stachel's detailed definition of the $(k,n)$-flexes \cite{stachel_proposal} introduced by Sabitov \cite{sabitov},
to resolve the dilemma of ending up with an infinite sequence of possible 
$(k,n)$-flexes, as pointed out by Stachel's example of the extended double-Watt mechanism \cite{stachel_aim}. 
Moreover, this algebraic approach also allows to take reality issues into account.

The presented computation of the redefined 
$(k,n)$-flexes used a resultant based projection method  \cite{melanova} to achieve the minimal parametrization
of the branches of algebraic space curves.  
For an efficient computation of the $(k,n)$-flexes we plan to resort for the underlying Puiseux expansion to the powerful mean of tropical geometry \cite{maclagan}, which was demonstrated in \cite{jensen} and also used in the kinematic context by Nayak \cite{nayak} for analyzing the configuration space of mechanisms. This approach is dedicated to future research as well as the study of the already mentioned  $p$-parametric flexes with $p>1$.

\paragraph{{\bf Acknowledgments}}
The author is supported by Grant No.\ F\,77 (SFB {\it Advanced Computational Design}, SP7) of the Austrian Science Fund FWF.

\newpage

\section*{Appendix A: Completing the discussion of Example 3}

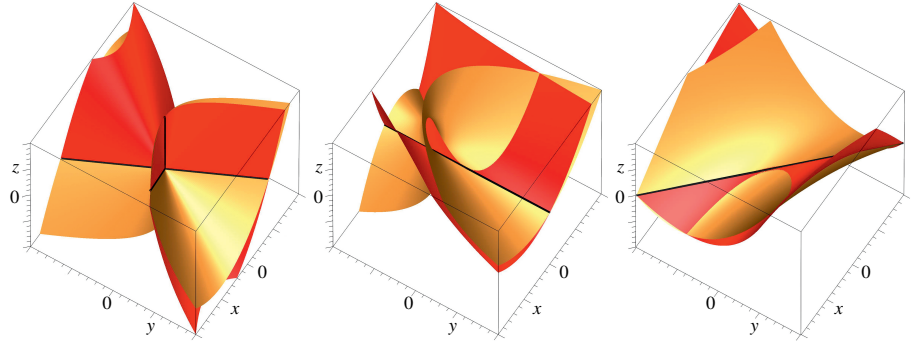
\begin{figure}[t]
\begin{center}
\begin{overpic}
    [height=45mm]{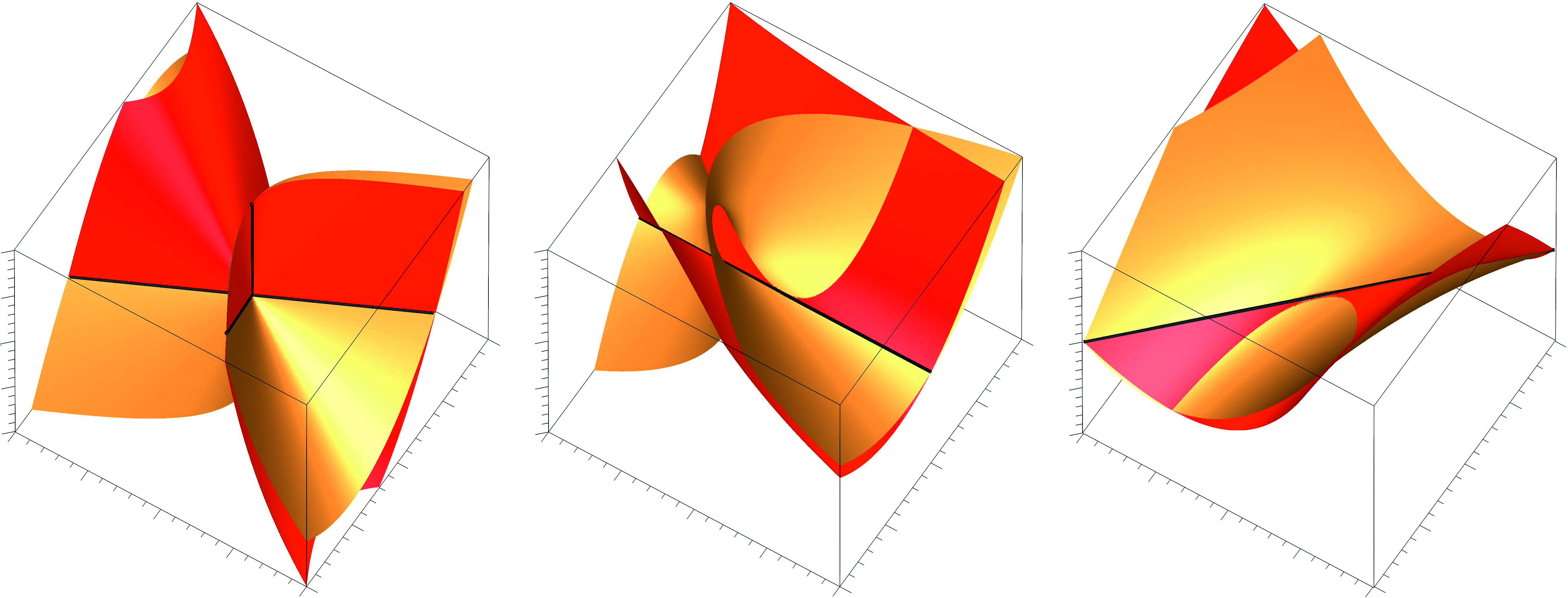}
\begin{scriptsize}
\put(-1.4,15.5){$0$}
\put(-1.2,19){$z$}
\put(32.8,15.5){$0$}
\put(33,19){$z$}
\put(66.8,15.5){$0$}
\put(67,19){$z$}
\put(14,1){$y$}
\put(9,3.5){$0$}
\put(48,1){$y$}
\put(43,3.5){$0$}
\put(82,1){$y$}
\put(77,3.5){$0$}
\put(23.5,3){$x$}
\put(26.5,7){$0$}
\put(57.5,3){$x$}
\put(60.5,7){$0$}
\put(91.5,3){$x$}
\put(94.5,7){$0$}
\end{scriptsize}     
  \end{overpic} 
\end{center}	
\caption{The quadric which corresponds to the first/second generator of $I_1$ is illustrated in red/yellow. Left: Case (1a), which equals case (2a). Center: Case (1b). Right: Case (2b).}
  \label{fig2}
\end{figure}

{\bf A) Discussion of $I_1$} \newline
Here we complete the discussion of the ideal $I_1$ started in Example \ref{ex3} by proceeding with the 
special cases excluded by  Eq.\ (\ref{eq:assume}).
\begin{enumerate}
\item 
$\lambda_1 = 3\lambda_2 - 1$: We proceed similar to the general case discussed in  Example \ref{ex3}. We end up with the following minimal parametrization for $P(x,y)$:
\begin{equation}
x(t)=t^3, \quad
y(t)=
    \left(\tfrac{4\lambda_2 - 2}{\lambda_2 - 1}\right)^{\tfrac{2}{3}}\tfrac{(\lambda_2 - 1)}{4\lambda_2 - 2}t^2+\ldots
\end{equation}
under the assumption 
\begin{equation}\label{eq:ass1}
    (2\lambda_2-1)(\lambda_2-1)\neq 0.
\end{equation}
Under this assumption we can also compute the minimal parametrization of $P(x,z)$ which reads as
\begin{equation}
x(t)=t^3, \quad z(t)=
\tfrac{1}{2}\left(\tfrac{4\lambda_2 - 2}{\lambda_2 - 1}\right)^{\tfrac{1}{3}} t^4+\ldots .
\end{equation}
If we put the resulting minimal parametrization of the space curve into $c_3$ we get a polynomial where $t^4$ factors out. Thus this branch implies a $(2,3)$-flex. Now we have to deal with the very special cases excluded by Eq.\ (\ref{eq:ass1}):
    \begin{enumerate}
        \item 
        $\lambda_2=1$: In this case we eliminate  the variables $x$ and $y$ by means of resultant, which yields
        \begin{equation}
            P(y,z):=(y^2 - 12yz + 4z^2)y^2 \quad
            P(x,z):=(x^2+4 z x-4 z^2) x^2
        \end{equation}
        in order to avoid the problem mentioned in Footnote \ref{fn1}. 
        It turns out that the space curve splits up into four straight lines, which can be parametrized as follows:
        \begin{equation}
        x_\mp(t) =-(2\mp 2\sqrt{2})t,\quad
        y_\mp(t)=(6\mp4\sqrt{2})t, \quad z(t)=t    
        \end{equation}
        and two times the $z$-axis given by $x(t)=y(t)=0$, $z(t)=t$. None of these four lines, which are illustrated in Fig.\ \ref{fig2}-left is tangent to $c_3$, thus the intersection multiplicity is always $1$. Thus no flex is implied by these four lines, as $n\geq k$ of Eq.\ (\ref{eq:flexk}) is violated.
        \item 
        $\lambda_2=\tfrac{1}{2}$: In this case we eliminate the variables $x$ and $z$ by means of resultant, which yields
        \begin{equation}
            P(x,y):=(2xy - 2y^2 - x)x \quad
            P(y,z):=(2y^3 - 4y^2z + 4yz - z)z
        \end{equation}
        because otherwise we would again run into the problem mentioned in Footnote \ref{fn1}. Now the space curve splits up into the $y$-axis ($x(t)=z(t)=0$, $y(t)=t$) and in a cubic curve touching the $y$-axis in the origin (cf.\ Fig.\ \ref{fig2}-center). The cubic has the minimal parametrization:
        \begin{equation}
         x(t)=-2t^2-4t^3+\ldots , \quad y(t)=t, \quad 
         z(t)=2t^3 +\ldots .
        \end{equation}
        Both linear branches intersect $c_3$ with multiplicity 2 implying two $(1,1)$-flexes. 
    \end{enumerate}
    \item
    $\lambda_1=\frac{1}{3}(8\lambda_2-2)$: In this case we compute again the projection $P(x,y)$ and $P(x,z)$ and we get finally the minimal parametrization of the space curve under the assumption
    \begin{equation}\label{eq:ass2}
        \lambda_2(\lambda_2-1)\neq 0
    \end{equation}
    as
    \begin{equation}
    x(t)=t^2, \quad 
    y(t)=-t^2+
    w_1t^3+\ldots ,\quad
    z(t)=-\tfrac{3\lambda_2}{\lambda_2 - 1}t^4+w_1t^5+\ldots 
    \end{equation}
    with
    $w_1:=\tfrac{2\lambda_2}{\lambda_2 - 1}\sqrt{\tfrac{3\lambda_2 - 3}{\lambda_2}}$. 
    Plugging this into $c_3$ shows again that this branch implies a $(2,3)$-flex. We remain with the discussion of the two very special cases excluded by Eq.\ (\ref{eq:ass2}):
    \begin{enumerate}
        \item 
        $\lambda_2=1$: This implies $\lambda_1=2$ and therefore we end up with exactly the same case as already discussed in case (1a). 
        \item 
        $\lambda_2=0$: In this case we compute the following two  projections: 
        \begin{equation}
            P(x,y):=(x + y)( x + y - 4y^2), \quad
            P(x,z):=z(4x^3 + 12xz - 16z^2 + z). 
        \end{equation}
        Now the space curve splits again up into a line ($x(t)=t,y(t)=-t,z(t)=0$) and a cubic curve, which touches the line in the origin  (cf.\ Fig.\ \ref{fig2}-right). The cubic has the minimal parametrization
        \begin{equation}
        x(t)=t, \quad 
        y(t)=-t +4t^2-32t^3+\ldots,\quad 
        z(t)=-4t^3+\ldots .
        \end{equation}
        Both linear branches intersect $c_3$ with multiplicity 2 which implies two $(1,1)$-flexes. 
        This completes the discussion of the ideal $I_1$.
    \end{enumerate}    
\end{enumerate}

\noindent
{\bf B) Discussion of $I_2$} \newline
Under the assumption of 
\begin{equation}\label{eq:ass3}
    (3\mu_1 + 8)(\mu_1 + 3)\neq 0
\end{equation}
we can compute the minimal parametrization of the two linear 
branches of the space curve using again the projections $P(x,y)$ and $P(x,z)$, which yields:
\begin{equation}
\begin{split}
&x(t)=t, \quad
y_\mp(t)=
-\tfrac{\mu_1+3\mp w_2}{\mu_1 + 3}t
-
\tfrac{(2\mu_1 + 5)(\mu_1 + 2)(6\mu_1 \mp 3w_2 + 16)}{(3\mu_1 + 8)(\mu_1 + 3)^2}t^2+\ldots, \\
&z_\mp(t)=\tfrac{3\mu_1 + 7 \mp w_2}{\mu_1 + 3}t^2+\ldots \quad \text{with} \quad
w_2:=\sqrt{(3\mu_1 + 8)(\mu_1 + 3)}.
\end{split}
\end{equation}
 Both linear branches intersect $c_2$ with multiplicity 2 which implies two $(1,1)$-flexes. 
 We proceed with the discussion of the two special cases excluded by Eq.\ (\ref{eq:ass3}):
 \begin{enumerate}
     \item 
     $\mu_1=-3$: In this case we can proceed as in the general case but we end up with one branch having the minimal parametrization
     \begin{equation}
         x(t)=t^3, \quad y(t)=\tfrac{1}{2}2^{\tfrac{1}{3}}t^2-\tfrac{1}{2}t^3+\tfrac{5}{12}2^{\tfrac{2}{3}}t^4+\ldots, \quad
         z(t)=\tfrac{1}{2}2^{\tfrac{2}{3}}t^4+\ldots
     \end{equation}
     Plugging this into $c_2$ shows that this branch implies a $(2,3)$-flex.
     \item 
      $\mu_1=-\tfrac{8}{3}$: Exactly the same holds as in the special case before. We only get the following minimal parametrization of the cuspidal branch: 
      \begin{equation}
          x(t)=t^2, \quad
          y(t)=t^2+2\sqrt{3}t^3-4t^4+\ldots, \quad
          z(t)=-3t^4+\ldots .
      \end{equation}
      which also implies a $(2,3)$-flex. 
      This completes the discussion of the ideal $I_2$.
 \end{enumerate}
\noindent
{\bf C) Discussion of $I_3$} \newline
This case corresponds to the removal of the first condition $c_1$ already briefly discussed in Example \ref{ex2}. Therefore, we expect to end up with two linear branches, which is also the case. By using the projections  $P(x,y)$ and $P(x,z)$ we get their minimal parametrizations as:
\begin{equation}
    x(t)=t,\quad
    y_\mp(t)=-(1\mp\sqrt{3})t-(4\mp2\sqrt{3})t^2+\ldots,\quad
    z_\mp(t)=(3\mp\sqrt{3})t^2+\ldots .
\end{equation}
Both linear branches intersect $c_1$ with multiplicity 2 which implies two $(1,1)$-flexes. This closes the discussion of all cases.


\section*{Appendix B: Immobile 4-bar mechanism}

We use again the projection approach of Melanova \cite{melanova} for the computation of the minimal parametrization of the branches of the algebraic curve through the origin, where we project to all coordinate planes containing the $b$-axis.

As equation $c_1=0$ only depends on the variables $a$ and $b$, it already corresponds to the projection of the space curve to the $ab$-plane; i.e.\ $P(a,b)=c_1$. Its minimal parametrization is given by 
\begin{equation}\label{eq:ab}
a(t)=-\tfrac{1}{2}t^2   -\tfrac{1}{8}t^4  + \ldots,\quad
b(t)=t.
\end{equation}
For the projection of the space curve to the $bc$-plane we first eliminate $a$ by means of resultant from the equations $c_1=c_2=0$. From the resulting expression $Res(c_1,c_2;a)$ and $c_3$ we eliminate $d$ again by means of resultant yielding:
\begin{equation}
 \begin{split}
  P(b,c):= &144b^4c^2 + 144b^2c^4 + 384b^4c + 456b^2c^3 + 256b^4 + 1376b^2c^2 +  \\
    &1225c^4 + 1024b^2c + 2800c^3 + 1600c^2.
 \end{split}    
\end{equation}
Its minimal parametrization splits up into the following two conjugate complex branches:
\begin{equation}\label{eq:bc}
b(t)=t,\quad
c(t)=\tfrac{-8 \pm 6I}{25}t^2 + \tfrac{-79 \pm 3I}{1250}t^4 + \ldots .
\end{equation}
Finally, we project the space curve to the $bd$-plane. This can be achieved by eliminating $c$ from $Res(c_1,c_2;a)$ and $c_3$ by means of resultant, which implies:
\begin{equation}
 \begin{split}
  P(b,d):= 
  &144b^4d^2 - 288b^3d^3 + 144b^2d^4 + 1024b^4 - 1408b^3d + 
  1968b^2d^2 - \\ &2080bd^3 + 
  1225d^4 + 11520b^2 - 5760bd + 7200d^2.
 \end{split}    
\end{equation}
Its minimal parametrization splits again up into two conjugate complex branches:
\begin{equation}\label{eq:bd}
b(t)=t,\quad
d(t)=\tfrac{2\pm 6I}{5}t+ \tfrac{6\pm 33I}{250}t^3 +\ldots .
\end{equation}
Then Eqs.\ (\ref{eq:ab}, \ref{eq:bc}, \ref{eq:bd}) imply the minimal parametrization of the two conjugate complex linear branches of the space curve given in Example \ref{ex4}. 
\end{document}